\documentclass{amsart}

\usepackage{tikz, hyperref, aliascnt, color, mathtools}
\definecolor{shade1}{rgb}{.2, .2, .2}
\definecolor{shade2}{rgb}{.5, .5, .5}
\definecolor{shade3}{rgb}{.8, .8, .8}
\definecolor{shade4}{rgb}{.95, .95, .95}
\usetikzlibrary{decorations.pathreplacing}
\usetikzlibrary{patterns}
\theoremstyle{plain}
\newtheorem{theorem}{Theorem}
\newtheorem*{theorem*}{Theorem}

\newtheorem{lemma}{Lemma}
\newaliascnt{corollary}{lemma}
\newtheorem{corollary}[corollary]{Corollary}
\aliascntresetthe{corollary}
\theoremstyle{definition}
\newaliascnt{definition}{lemma}
\newtheorem{definition}[definition]{Definition}
\aliascntresetthe{definition}

\def\sek~{\S{}}

\def\H{\mathcal{H}}

\title[Residual generic ergodicity]{Residual generic ergodicity of periodic group extensions over translation surfaces}

\author{David Ralston}
\address{SUNY College at Old Westbury, Mathematics/CIS Department, P.O. Box 210, Old Westbury, NY 11568, USA}
\email{ralstond@oldwestbury.edu}

\author{Serge Troubetzkoy}
\address{Aix Marseille Universit\'e, CNRS, Centrale Marseille, I2M, UMR
  7373, 13453 Marseille, France}
  
\curraddr{ I2M, Luminy, Case 907, F-13288 Marseille CEDEX 9, France}
 \email{serge.troubetzkoy@univ-amu.fr}

\thanks{This work was partially supported by ANR Perturbations}
\date{\today}

\begin{document}

\begin{abstract}
Continuing the work in \cite{ergodic-infinite}, we show that within each stratum of translation surfaces, there is a residual set of surfaces for which the geodesic flow in almost every direction is ergodic for almost-every periodic group extension produced using a technique referred to as \emph{cuts}.
\end{abstract}
\maketitle

\section{Introduction}

A translation surface is a real 2-dimensional manifold with conical singularities equipped with an atlas for which transition functions are translations.  The simplest example of a compact translation surface is the 2 dimensional torus.  It is a classical result that the straight-line flow on the torus satisfies a dichotomy: for rational directions it is completely periodic, while for irrational directions it is uniquely ergodic. There has been considerable work in extending and generalizing this result to other translation surfaces; here we list the major achievements.  Veech has shown that this dichotomy hold for any compact translation surface whose group of affine diffeormorphism is a lattice in $SL(2,\mathbb{R})$ \cite{V}. Such surfaces  are now known as {\em Veech surfaces}, this class includes all regular polygons with an even number of sides, with  parallel sides identified, as well as compact square-tiled translation surfaces and other examples. Kerckhoff, Masur, and Smillie have shown that on any compact translation surface, the geodesic flow in almost every direction is uniquely ergodic \cite{1986}. Avila and Forni have shown that for almost every compact translation surface of genus at least 2, for almost every direction the geodesic flow is in fact weak-mixing \cite{AF}. 

On the other hand, much less is known for noncompact translation surfaces.  For aperiodic infinite translation surfaces, recurrence is known for generic translation surfaces in certain classes \cite{TTT,TT}. A bit more is known for periodic infinite translation surfaces.  Hubert, Hooper, and Weiss constructed an example of an infinite translation surface for which the geodesic flow is ergodic
in almost every direction \cite{HHW}.  One of the most famous infinite translation surfaces arises from the classical Ehrenfest wind-tree model, introduced by Paul and Tatiana Ehrenfest in 1912 \cite{EE}.  Hubert, Leli\`evre, and Troubetzkoy, and Avila and Hubert proved the recurrence of the geodesic flow for periodic wind-tree model  \cite{HLT,AH} and the diffusion rate was calculated
by Delecroix, Hubert, and Leli\`evre to be 2/3 \cite{DHL}. Fra{\c{c}}zek and Ulcigrai have shown that the geodesic flow for wind-tree model is not ergodic for almost all directions \cite{FU}. In \cite{ergodic-infinite} we showed that for square tiled translation surfaces which have a single cylinder direction,
almost every infinite group extensions of geodesic flows produced by cuts is ergodic in almost every direction.

In this article  we extend this result to the following: there is a residual set of surfaces for which the geodesic flow in almost every direction is ergodic for generic periodic group extensions produced by cuts. The techniques used are the same ones as in \cite{ergodic-infinite}: essential values and approximation by purely periodic directions.

\section{Background and Details of the Construction}\label{section - background}
A translation surface is a Riemannian surface $M$, a finite or countable set of points $
\Sigma := \{ x_1,x_2,\dots \} \subset M$  (the singularity set) and an open cover of $M \setminus \Sigma$ by sets $\{ U_{\alpha} \}$, together with charts $\phi_{\alpha} : U_{\alpha} \to \mathbb{R}^2$ such that for all $\alpha,\beta$ with $U_\alpha \cap U_\beta \ne \emptyset$ the map $\phi_\alpha \circ \phi_\beta^{-1}$ is a translation.  The Euclidean metric and Lebesgue measure on $\mathbb{R}^2$ may be pulled back to $M$ via these charts except at 
singular points $x_i \in \Sigma$ with cone angle $2 \pi (d_i+1)$ for $d_i \in \mathbb{Z}^+$; this set $\Sigma$ is finite if $M$ is compact.  We denote by $\mu$ the pull-back of Lebesgue measure to $M$ (if $\mu(M)< \infty$, we normalize to that $\mu(M)=1$).  The collection of all compact translation surfaces with $n$ singularities of cone angles $2 \pi (d_i+1)$ for $i=1,2,\ldots,n$ (where the $d_i$ are monotone decreasing by convention) is referred to as a \emph{stratum}, and denoted by $\mathcal{H}(d_1,\ldots,d_n)$.  Any compact translation surface $M$ may be viewed (in a non unique way) as a polygon (not necessarily convex) in $\mathbb{R}^2$, where pairs of parallel congruent sides are identified (see e.g. \cite[\S12]{Viana2006}).
 The precise parameterization of strata that we will use is that of the space of {\em zippered rectangles}, which we will define later.

 When $M$ is a $k$-fold cover of the torus ramified over a point, then $M$ is called \emph{square-tiled}.  On each flat surface $M$ we also choose a marked point $p$ in a specified chart $U_{\alpha}$, and we require only $p \notin \Sigma$.

On $M$, with the flat metric lifted from the charts, geodesics that do not go through singularities project to straight lines in the plane in a fixed direction, and such geodesics are either periodic or simple (do not intersect themselves).  If the geodesic in direction $\theta$ starting at $x \in M$ does not intersect with $\Sigma$ then we define the geodesic flow $\varphi_t^{\theta}(x)$ to the point obtained after moving in the direction $\theta$ for a time $t$ starting at $x$.  If $\theta$ is fixed,  we may simply write $\varphi_t(x)$.  If the geodesic hits a singularity, then the geodesic flow $\varphi_t(x)$ is defined up to this time. The geodesic flow preserves the surface area measure $\mu$.  For a survey on translation surfaces, see \cite{Viana2006, MR2261104}.

We may lift generic points and geodesics from the plane to $M$;  more precisely consider the marked point $p \in U_{\alpha}$ and  $p' :=  \phi_{\alpha}(p) \in \mathbb{R}^2$
Let $x \in \mathbb{R}^2$ be arbitrary, and consider the line segment connecting $p'$ to $x$. Let $\ell$ denote the length of this segment, note that part of this  segment is in $\phi_{\alpha}(U_{\alpha})$.  Lift this part  to $U_{\alpha}$ and, if possible (i.e. if the procedure does not run into a singularity),
continue this line in $M$ until it is of length $\ell$, in which case 
the line segment from $p'$ to $x$ on $\mathbb{R}^2$ {\em lifts} to a unique geodesic on $M$. This lift exists for almost
all $x \in \mathbb{R}^2$.  

 \begin{definition}
A \emph{collection of cuts} $\gamma$ is a finite collection of pairs of parallel congruent line segments in the plane.  Specifically, let $N \in \mathbb{Z}^+$, and for each $i=1,2,\ldots,N$, let
\[ \left( \textbf{x}_i^+, \textbf{x}_i^- \right) \coloneqq \left( (x_i^+, y_i^+), (x_i^-, y_i^-) \right) \in \mathbb{R}^2 \times \mathbb{R}^2\] be a pair of points in the plane, $\zeta_i \in \mathbb{S}^1$, and $\ell_i \in \mathbb{R}^+$.  Denote by $\gamma_i^+$ and $\gamma_i^-$ the pair of geodesics beginning at the two points $\mathbf{x}^{\pm}_i$, in direction $\zeta_i$ and of length $\ell_i$.  All possible finite collections of cuts are parameterized by the $\sigma$-finite measure space
\begin{equation}\label{eqn - space of cuts} \Gamma \coloneqq \bigcup_{N=1}^{\infty} \left(  \mathbb{R}^2 \times \mathbb{R}^2 \times \mathbb{S}^1 \times \mathbb{R}^+\right)^N ,
\end{equation}
equipped with Lebesgue measure.  We will generally not distinguish between the points $\textbf{x}_i^{\pm}$ in the plane and their lifts to a surface $M$, and similarly for the segments $\gamma_i^{\pm}$, allowing us to more easily speak of placing the same collection of cuts on two different surfaces.
\end{definition}

Given a compact translation surface $M$, we will form \textit{skew products} over $M$ in the following manner: let $G$ be a locally compact Abelian group with translation-invariant metric $\| \cdot \|$, whose identity element we denote $e_G$, and let $\gamma$ be a collection of $N$ cuts.  Lift each segment in $\gamma$ to $M$; if some $\gamma_i^{\pm}$ encounters a singularity on $M$ then the lift will not be unique, but there will be only finitely many continuations, and for the time being we may simply make a choice of continuation.  Let $\{g_1,g_2,\ldots,g_N\}$ be a finite set in $G$.  Define $f:M \rightarrow G$ to take value $g_i$ on the lift of each $\gamma_i^{+}$, value $-g_i$ on the lift of $\gamma_i^-$, and value $e_G$ elsewhere.  Wherever two (or more) lifts of $\gamma_i$ cross, define the function at these points to be the sum of the function on those segments.  Given a point $x \in M$ and a direction $\theta$ which is not equal to any $\zeta_i$, denote
\begin{equation}\label{eqn - sums} S_t(x) \coloneqq S_t(x,M,\theta) \coloneqq \sum_{\substack{ 0 \leq s <t \\ f(\varphi_s^{\theta}(x)) \neq e_G}}f(\varphi^{\theta}_s(x)),\end{equation}
and with $M$ and $\theta$ fixed,
\begin{equation}\label{eqn - balls for essential values}
S_t^{-1} \left( B_{\epsilon}(g) \right) \coloneqq \left\{ x \in M : \|S_t(x) - g \| < \epsilon \right\}.
\end{equation}

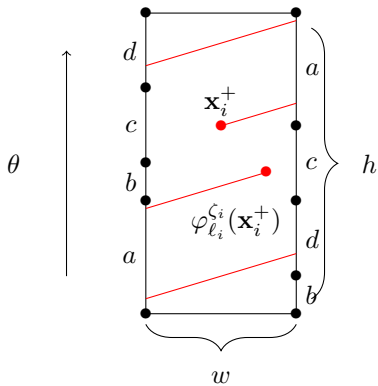
\begin{figure}
\center{\begin{tikzpicture}
\draw (0,0)-- node[midway, left]{$a$}(0,1.5)-- node[midway, left]{$b$}(0,2)-- node[midway, left]{$c$}(0,3)-- node[midway, left]{$d$}(0,4)--(2,4)-- node[midway, right]{$a$}(2,2.5)-- node[midway,right]{$c$}(2,1.5)--node[midway,right]{$d$}(2,.5)-- node[midway,right]{$b$}(2,0)--(0,0);

\draw [color=red] (1,2.5) -- (2,2.8);
\draw [color=red] (0,.2) -- (2,.8);
\draw [color=red] (0,3.3) -- (2,3.9);
\draw [color=red] (0,1.4) -- (1.6,1.88);
\node at (1,2.5)[above] {$\textbf{x}_i^+$};
\node at (1,2.5)[color=red] {\textbullet};
\node at (1.2,1.6)[below] {$\varphi_{\ell_i}^{\zeta_i}(\textbf{x}_i^+)$};
\node at (1.6,1.88)[color=red] {\textbullet};
\foreach \Point in {(0,0), (0,1.5), (0,2), (0,3), (2,.5), (2,1.5), (2,2.5), (2,0), (0,4), (2,4)}
    \node at \Point {\textbullet};
    
    \draw [decorate, decoration={brace,amplitude=10pt, mirror},xshift=0pt,yshift=0pt] (2.2,.2)--  (2.2,3.8) node[midway, xshift=22pt] {$h$};
\draw [decorate, decoration={brace,amplitude=10pt, mirror}, yshift=-4pt] (0,0) -- (2,0)node[midway,yshift=-20pt] {$w$};
\draw [xshift=-30pt](0,.5) edge[->] node[midway, xshift=-20pt]{$\theta$} (0,3.5) ;

\end{tikzpicture}}
\caption{\label{figure - single cylinder}The single-cylinder representation for a surface $M$ with completely periodic flow in direction $\theta$ with period $h$, with the flow represented as vertical.  The length $w=\mu(M)/h$ is called the width of the cylinder.  Parallel saddle connections of equal lengths are identified by the letters $a$ through $d$, and the top and bottom are identified. A sample cut on the surface is the lift of $\gamma_i^+$, beginning at $\textbf{x}_i^+$ (the lift of $(x_i^+,y_i^+)$ to $M$), in direction $\zeta_i$ and of length $\ell_i$.  A corresponding cut, parallel and of equal length, is the lift of $\gamma_i^-$.}
\end{figure}

Then the transformation $\tilde{\varphi}_t$, from $\tilde{M}=M \times G$ to itself, is given by
\[\tilde{\varphi}_t(x,g) = (\varphi_t(x),g+S_t(x)),\] and it is this continuous time flow which we study; $\tilde{\varphi}_t$ preserves the measure $\tilde{\mu}=\mu \times \nu$, where $\nu$ is the Haar measure on $G$.  See \autoref{figure - single cylinder} to see how one single cut $\gamma_i^+$ might be appear on a surface with a \emph{single-cylinder direction}, i.e. a direction $\theta$ so that 
the flow on $M$ in this direction is completely periodic, with all points having the same period $h$, such that there is an open, connected set which is invariant under the flow in this direction, contains no singularities, and is of full measure.
  The surface $M$ may thus be represented as a strip with the elements of $\Sigma$ restricted to the sides parallel to the flow.  Geodesics which connect two singularities are called \emph{saddle connections}.

In order for the flow on $\tilde{M}$ to be ergodic in direction $\theta$, it is necessary that the projection of the cuts to a segment in a direction transverse to $\theta$ has an average value of $e_G$ (for example, projected to the base of a single-cylinder).  Note that by assigning the values $\pm g_i$ on parallel segments of equal length, we have ensured that this condition is always met.

\begin{definition}[Generically ergodic extensions]\label{definition - generic extensions}
Suppose that $M$ is a compact translation surface so that for any locally compact Abelian group $G$, for almost every collection of cuts and almost every $\theta \in \mathbb{S}^1$, the flow $\tilde{\varphi}_t^{\theta}$ on $\tilde{M}$ is ergodic for any choice of values $\{g_1,g_2,\ldots,g_N\}$ on the cuts which generate a dense subgroup of $G$.  Then $M$ will be said to \emph{produce generically ergodic extensions via cuts}. Let $\mathcal{E}$  be the set of translation surfaces $M$ which produce generically ergodic extensions via cuts.
\end{definition}

One central result of \cite{ergodic-infinite} is that if $M$ is a square-tiled surface with a single-cylinder direction, then $M \in \mathcal{E}$.  It is of interest to know how typical this property is for translation surfaces in general.  In this article we are not concerned with answering this question regarding measure-theoretic generic $M$, but rather we answer the question of \emph{topological} genericity.  In a Baire space $Y$, a set $A \subseteq Y$ is called $G_{\delta}$ if it may be written as the countable intersection of open, dense sets, and $A$ is called \emph{residual} if it contains a $G_{\delta}$ set. 

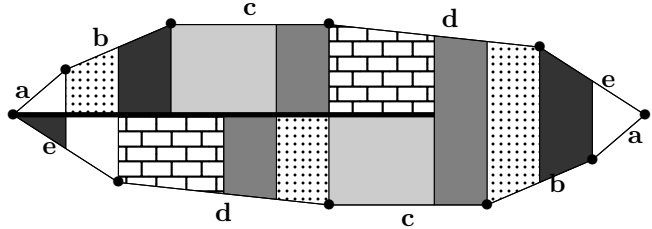
\begin{figure}[tbh]
\center{\begin{tikzpicture}[xscale = .7, yscale=.15]
\draw [draw=black, fill=white] (0,0) --(1,4) -- (1,0)--(0,0);
\draw [draw=black, fill=white] (1,-3)--(1,0)--(2,0)--(2,-6)--(1,-3);
\draw [draw=black, fill=white] (11,-4)--(11,3)--(12,0)--(11,-4);
\draw [draw=black, pattern=dots] (1,0)--(2,0)--(2,6)--(1,4)--(1,0);
\draw [draw=black, pattern=dots] (5,-7.5)--(5,0)--(6,0)--(6,-8)--(5,-7.5);
\draw [draw=black, pattern=dots] (9,6.5)--(10,6)--(10,-6)--(9,-8)--(9,6.5);
\draw [draw=black, fill=shade1] (1,-3)--(0,0)--(1,0)--(1,-3);
\draw [draw=black, fill=shade1] (2,0)--(2,6)--(3,8)--(3,0)--(2,0);
\draw [draw=black, fill=shade1] (10,-6)--(10,6)--(11,3)--(11,-4)--(10,-6);
\draw [draw=black, fill=shade3] (3,0)--(3,8)--(5,8)--(5,0)--(3,0);
\draw [draw=black, fill=shade3] (6,-8)--(6,0)--(8,0)--(8,-8)--(6,-8);
\draw [draw=black, fill=shade2] (5,0)--(5,8)--(6,8)--(6,0)--(5,0);
\draw [draw=black, fill=shade2] (8,-8)--(8,7)--(9,6.5)--(9,-8)--(8,-8);
\draw [draw=black, fill=shade2] (4,-7)--(4,0)--(5,0)--(5,-7.5)--(4,-7);
\draw [draw=black, pattern=bricks] (6,0)--(8,0)--(8,7)--(6,8)--(6,0);
\draw [draw=black, pattern=bricks] (2,0)--(4,0)--(4,-7)--(2,-6)--(2,0);
\foreach \Point in {(0,0)),(1,4),(3,8),(6,8),(10,6),(12,0),(11,-4),(9,-8),(6,-8),(2,-6)}
    \node at \Point {\textbullet};
\draw (0,0)--node[midway, left]{\textbf{a}}(1,4)--node[midway, left,yshift=4pt]{\textbf{b}}(3,8)--node[midway, above]{\textbf{c}}(6,8)--node[midway, above, xshift=6pt]{\textbf{d}}(10,6)--node[midway, right]{\textbf{e}}(12,0)--node[midway, right]{\textbf{a}}(11,-4)--node[midway, right]{\textbf{b}}(9,-8)--node[midway, below]{\textbf{c}}(6,-8)--node[midway, below]{\textbf{d}}(2,-6)--node[midway, left]{\textbf{e}}(0,0);
\draw [line width=2] (0,0)--(8,0);
\end{tikzpicture}}
\caption{\label{figure - polygon} A polygon (e.g. a decagon) with pairs of parallel sides defines a translation structure (side identifications are marked).  By defining a segment transverse to the vertical flow, we define a first-return map to this segment.
}
\end{figure}

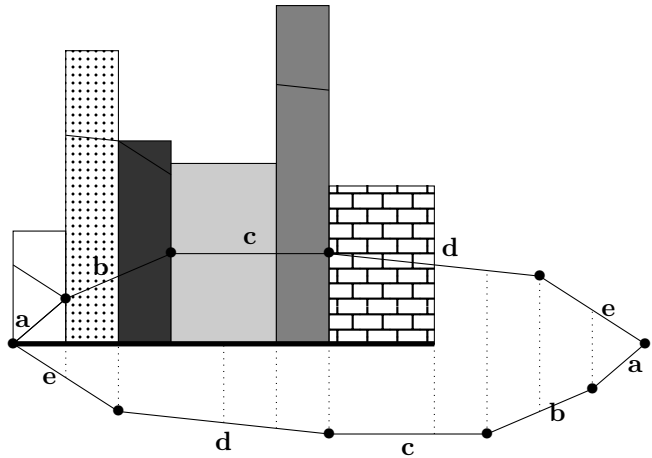
\begin{figure}[tbh]
\center{\begin{tikzpicture}[xscale = .7, yscale=.15]
\draw [draw=black, fill =white](0,0)--(0,10)--(1,10)--(1,0)--(0,0);
\draw [draw=black, pattern = dots](1,0)--(2,0)--(2,26)--(1,26)--(1,0);
\draw [draw=black, fill=shade1](2,0)--(3,0) --(3,18)--(2,18)--(2,0);
\draw [draw=black,fill=shade3](3,0)--(5,0)--(5,16)--(3,16)--(3,0);
\draw [draw=black,fill=shade2] (5,0)--(6,0)--(6,30)--(5,30)--(5,0);
\draw [draw=black,pattern=bricks] (6,0)--(8,0)--(8,14)--(6,14)--(6,0);
\foreach \Point in {(0,0)),(1,4),(3,8),(6,8),(10,6),(12,0),(11,-4),(9,-8),(6,-8),(2,-6)}
    \node at \Point {\textbullet};
\draw (0,0)--node[midway, left]{\textbf{a}}(1,4)--node[midway, left,yshift=4pt]{\textbf{b}}(3,8)--node[midway, above]{\textbf{c}}(6,8)--node[midway, above, xshift=6pt]{\textbf{d}}(10,6)--node[midway, right]{\textbf{e}}(12,0)--node[midway, right]{\textbf{a}}(11,-4)--node[midway, right]{\textbf{b}}(9,-8)--node[midway, below]{\textbf{c}}(6,-8)--node[midway, below]{\textbf{d}}(2,-6)--node[midway, left]{\textbf{e}}(0,0);
\draw [line width=2] (0,0)--(8,0);
\draw (0,0)--(1,4)--(0,7);
\draw (1,18.5)--(2,18);
\draw (2,18)--(3,15);
\draw (5,23)--(6,22.5);
\draw [dotted](1,-3)--(1,0);
\draw [dotted](2,-6)--(2,0);
\draw [dotted](4,-7)--(4,0);
\draw [dotted](5,-7.5)--(5,0);
\draw [dotted](6,-8)--(6,0);
\draw [dotted](8,-8)--(8,0);
\draw [dotted](9,-8)--(9,6.5);
\draw [dotted](10,-6)--(10,6);
\draw [dotted] (11,-4)--(11,3);
\end{tikzpicture}}
\caption{\label{figure - rectangles} The first-return flow to the segment in \autoref{figure - polygon} may be rearranged as a collection of rectangles rather than an arbitrary polygon, where the tops of rectangles are identified to the bottom of the segment, inducing an \emph{interval exchange transformation}.  
}\end{figure}

We will work in the space of {\it zippered rectangles}; this space was introduced by Veech \cite{MR644019}, and we follow the presentation of Zorich \cite[\S5.7]{MR2261104}. 
Fix integers $d_i \ge 1$, and the corresponding stratum $\H(d_1,\dots,d_n)$.  Surfaces in this stratum have $n$ singular points, their  genus $g$ satisfies
$ \sum_{i=1}^n d_i = 2g-2$. 
Fix a translation surface $M$ in this stratum and a direction $\theta$ without saddle connections, this direction will be called
{\em vertical}, the perpendicular direction  will be called {\em horizontal}. For the rest of this article the term {\em translation surface} will correspond to
the pair $(M,\theta)$.
Suppose that $X$  is a horizontal section. The first return map to $X$ is always an interval exchange transformation. We suppose
furthermore that $X$ satisfies Convention 2 of \cite{MR2261104}: the first return of the vertical flow to $X$ has the minimal possible number
$$N := \sum_{}(d_j+1) + 1 = 2g + n -1$$
of subintervals under exchange.
Furthermore we will always place the left endpoint of $X$ at a singularity, the convention leaves a discrete choice for the position of the right endpoint.
The first return map decomposes $M$ into rectangles: these rectangles together with their identifications are referred to as {\em zippered rectangles}.  See \autoref{figure - polygon} for a convex polygon with vertical flow in which we induce the first-return flow to a horizontal segment, and \autoref{figure - rectangles} for the rearrangement of the figure as a collection of rectangles.

The coordinates of zippered rectangle start with the usual coordinates of interval exchange:
the widths $\lambda := (\lambda_1,\dots,\lambda_N)$ of the rectangles (intervals), and
the permutation $\pi \in S_N$ describing the order of the exchanged intervals.
Next we need the heights of $h := (h_1,\dots,h_N)$ the rectangles,  and $a := (a_1,\dots,a_N)$ the altitudes responsible for
the position of the singularities (note that the left endpoint of $X$ starts at a singularity, and the vertical direction does not have any saddle connections,
thus it gets mapped to $X$ without hitting another singularity). We zip neighboring rectangles $R_j$ and $R_{j+1}$ from the bottom
upto the height $a_j$ and then
split them. To complete the zipping, neighboring rectangles $R_{\pi^{-1}(j)}$ and $R_{\pi^{-1}(j+1)}$ are zipped from the top down to the height $h_{\pi^{-1}(j)} -a_{\pi^{-1}(j)}$.

The parameters $h_j$ and $a_j$ are not independent, they satisfy certain linear equations and linear inequalities (the picture
zips up completely, and the zippers do not exceed the top of the rectangles). Continuously varying the parameters
$\lambda,h,a$ while respecting the linear relations yields  coordinates in the stratum.  Veech showed that   given an interval exchange
transformation without saddle connections, one can always construct a translation surface $M$ and a horizontal segment $X \subset M$ 
such that that the first return map of the vertical flow to $X$ gives the initial interval exchange transformation. Moreover, there is a $n$-dimensional family of such flat surfaces, suspensions over the interval exchange transformation. The continuous variation  is a completion of the procedure, 
it includes directions with saddle connections, we will refer to the completion as the {\it boundary} of the space of zippered rectangles.

The choice of $X$ completely determines the
decomposition, thus one would like a canonical way to choose $X$. 
There is an almost canonical way to choose $X$, namely place the left endpoint of $X$ at one on the singularities and choose the
length of $X$ in such a way that $X$ is the shortest interval satisfying Convention 2 and the length of $X$ is at least 1. 
 We will not need to use it, but we remark that there is a constructive way to find $X$, start with a sufficiently long interval $X$ have its left end point at a singular point and
 satisfying Convention 2. Apply the ``slow'' Rauzy induction as long as the resulting interval has length at least one.  For almost all
 flat surface this procedure will stop after a finite number of steps, yielding the desired interval $X$.

Since $M$ has a finite number of singularities, and each singularity has a finite number of horizontal prongs, this yields  $\sum_{i=1}^n (d_i +1)= 2g-2+n$
choices for $X$.
Thus the space of zippered rectangles can be essentially viewed as a ramified covering of the corresponding connected component
$\H^{comp}(d_1,\dots,d_n)$ of the stratum.  Supposing that the area of the surfaces is one, defined by the condition $\lambda \cdot h = 1$ (a final
linear equation),
yields the space of zippered rectangles  of area one covering the space $\H^{comp}(d_1,\dots,d_n)$ of translation surfaces of area one.
Slightly abusing notation, we will refer to the space of zippered rectangles as  ${\H}(d_1,\dots,d_n)$.
Summarizing, we have shown that the space $\H(d_1,\ldots,d_n)$ is the product of a discrete countable parameter, the permutation (Rauzy class)
and the closure of an open subset of $\mathbb{R}^{3N}$ (given by the linear equations and linear inequalities), yielding 

\begin{lemma}\label{lemma - baire}
Each stratum $\H(d_1,\ldots,d_n)$ is a separable, metric Baire space.
\end{lemma}



\begin{figure}[tb]
\centering
\begin{minipage}{.4 \linewidth}
\center{\begin{tikzpicture}
\draw (0,0)-- node[midway, right]{$a$}(0,1.5)-- node[midway, right]{$b$}(0,2)-- node[midway, right]{$c$}(0,3)-- node[midway, right]{$d$}(0,4)--(2,4)-- node[midway, left]{$a$}(2,2.5)-- node[midway,left]{$c$}(2,1.5)--node[midway,left]{$d$}(2,.5)-- node[midway,left]{$b$}(2,0)--(0,0);
\draw [decorate, decoration={brace,amplitude=10pt, mirror},xshift=0pt,yshift=0pt] (2.2,.2)--  (2.2,3.8) node[midway, xshift=22pt] {$h$};
\draw [decorate, decoration={brace,amplitude=10pt, mirror}, yshift=-4pt] (0,0) -- (2,0)node[midway,yshift=-20pt] {$w$};
\draw [xshift=-30pt](0,.5) edge[->] node[midway, xshift=-20pt]{$\theta$} (0,3.5) ;
\foreach \Point in {(0,0), (0,1.5), (0,2), (0,3), (2,.5), (2,1.5), (2,2.5), (2,0), (0,4), (2,4)}
    \node at \Point {\textbullet};
\end{tikzpicture}}
\end{minipage}%
\hspace{.1 \linewidth}
\begin{minipage}{.4 \linewidth}
\centering
{\begin{tikzpicture}
\filldraw[color=shade2](0,0)--(-.1,1.6)--(0,1.9)--(0,3)--(.1,4.1)--(0,0);
\filldraw[color=shade2](2,4)--(2.1,2.4)--(2.1,1.3)--(2,.2)--(1.9,-.1)--(2,4);
\draw[dotted](0,0) -- (.1,4.1);
\draw[dotted](1.9,-.1) -- (2,4);
\draw (0,0)--(-.1,1.6)--(0,1.9)--(0,3)--(.1,4.1)--(2,4)--(2.1,2.4)--(2.1,1.3)--(2,.2)--(1.9,-.1)--(0,0);
\draw [decorate, decoration={brace,amplitude=10pt, mirror},xshift=4pt,yshift=0pt] (2.2,.2)-- (2.3,3.8) node[midway, xshift=22pt]{$h^*$};
\draw [decorate, decoration={brace,amplitude=10pt, mirror}, yshift=-4pt] (0,0) -- node[midway,yshift=-25pt]{$w^*$}(2,0) ;
\node at (1,2){$A$};
\draw [xshift=-20pt](0,0) edge[->] node[midway, xshift=-15pt]{$\theta^*$}(.1,4.1) ;
\foreach \Point in {(0,0), (-.1,1.6), (0,1.9), (0,3),(.1,4.1),(2,4), (2,.2), (2.1,1.3), (2.1,2.4), (1.9,-.1)}
    \node at \Point {\textbullet};
\end{tikzpicture}}
\end{minipage}
\caption{\label{figure - almost single cylinder}Left: a single-cylinder direction $\theta$ (not necessarily close to $\alpha_m$) on the surface $M_m$ with width $w$ and height $h$.  Parallel saddle connections of equal lengths are identified by the letters $a$ through $d$, and the top and bottom are identified.
Right: the vectors corresponding to the saddle connections  have been allowed to change slightly; such adjustments give a local metric within each stratum. The flow in direction $\theta^*$ (which is very close to $\theta$) will be periodic with nearly the same period $h^*$ on a set $A$ of large measure (excluding the shaded region).}
\end{figure}
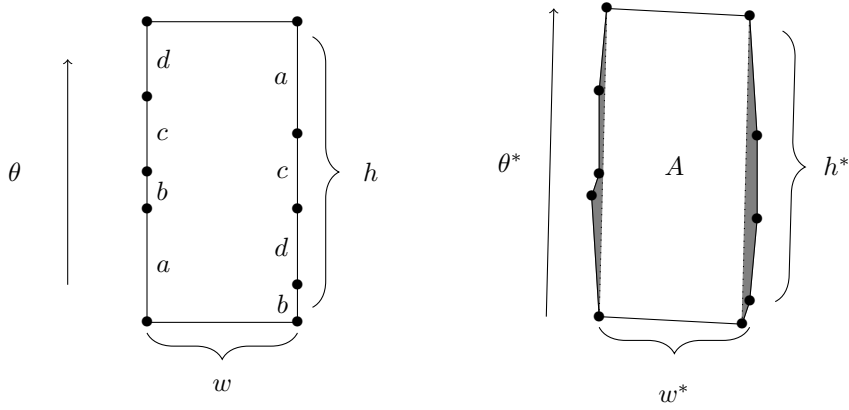

Note that if we have a translation surface $M$ and a direction $\theta$ which decomposes $M$ into a single cylinder, then this direction on this surface is in the boudary 
of the space of zippered rectangles since the direction of the flow contains saddle connections).  However, for \emph{nearby} directions $\theta'$, the pair $(M, \theta')$ will appear in the interior of the space of zippered rectangles (there are only countably many directions which admit saddle connections in any translation surface).  See \autoref{figure - almost single cylinder}, left.

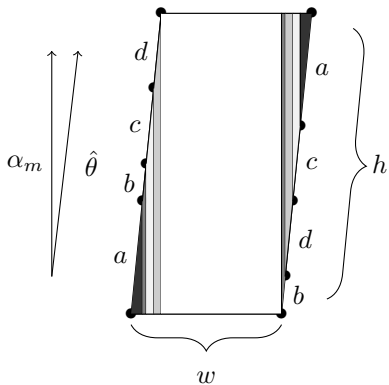
\begin{figure}[th]
\centering
{\begin{tikzpicture}
\draw (0,0)-- node[midway, left]{$a$}(.150,1.5)-- node[midway, left]{$b$}(.20,2)-- node[midway, left]{$c$}(.30,3)-- node[midway, left]{$d$}(.40,4)--(2.4,4)-- node[midway, right]{$a$}(2.25,2.5)-- node[midway,right]{$c$}(2.15,1.5)--node[midway,right]{$d$}(2.05,.5)-- node[midway,right]{$b$}(2,0)--(0,0);
\draw [decorate, decoration={brace,amplitude=10pt, mirror},xshift=20pt,yshift=0pt] (1.9,.2)--  (2.25,3.8) node[midway, xshift=15pt] {$h$};
\draw [decorate, decoration={brace,amplitude=10pt, mirror}, yshift=-4pt] (0,0) -- (2,0)node[midway,yshift=-20pt] {$w$};
\draw [xshift=-30pt](0,.5) edge[->] node[midway, xshift=10pt]{$\hat{\theta}$} (0.35,3.5) ;
\draw [xshift=-30pt] (0,.5) edge[->] node[midway, xshift=-10pt]{$\alpha_m$} (0,3.5);
\foreach \Point in {(0,0), (.150,1.5), (.20,2), (.30,3), (2.05,.5), (2.15,1.5), (2.25,2.5), (2,0), (.40,4), (2.4,4)}
    \node at \Point {\textbullet};
   \draw  [draw=black, fill=shade1] (0,0) -- (.15,0) -- (.15,1.5);
    \draw  [draw=black, fill=shade2]  (.15,0) -- (.2,0) -- (.2,2) -- (.15,1.5);
     \draw  [draw=black, fill=shade4]  (.2,0) -- (.3,0) -- (.30,3) -- (.2,2);
      \draw  [draw=black, fill=shade3]  (.3,0) -- (.4,0) -- (.4,4) -- (.3,3);
       \draw  [draw=black, fill=white]  (.4,0) -- (2,0) -- (2,4) -- (.4,4);
        \draw  [draw=black, fill=shade2]  (2,0) -- (2,4) -- (2.05,4) -- (2.05,.5) ;
      \draw  [draw=black, fill=shade3]   (2.05,.5) -- (2.15,1.5) -- (2.15,4) -- (2.05,4) ;
           \draw  [draw=black, fill=shade4]   (2.15,1.5) -- (2.25,2.5) -- (2.25,4) -- (2.15,4) ;
     \draw[draw=black, fill=shade1] (2.4,4)-- (2.25,4) -- (2.25,2.5);
\end{tikzpicture}}
\caption{\label{figure - dense collection} The zippered rectangle representation where the vertical direction $\alpha_m$ has no saddle connections and is close to the original single-cylinder direction $\hat{\theta}$.  Parallel saddle connections of equal lengths are identified by the letters $a$ through $d$, and the top and bottom are identified via an interval exchange.}
\end{figure}

\begin{lemma}\label{corollary - silly}
Within each stratum $\H(d_1, \ldots, d_n)$, there is a countable, dense set of $(M_m, \alpha_m)$, where $M_m$ is compact, square-tiled, and has a single-cylinder direction.
\begin{proof}
Lemma 18 and Remark 7 of  \cite{kontsevich-zorich} show that the set of compact surfaces with single-cylinder directions is dense in each stratum, where denseness is not in the set of zippered rectangles but in the sense of closeness of Abelian differentials. The constructed surfaces have the vertical direction as single-cylinder (which is forbidden in the space of zippered rectangles); we denote this direction as $\hat{\theta}$.

By making arbitrarily small changes to the vectors of the vertical saddle connections, we ensure that the lengths of all vertical saddle connections are rational multiples of the width of the cylinder without changing the fact that the vertical flow is a single-cylinder direction.  By rescaling the width by an arbitrarily small amount, we may ensure all such lengths are in fact rational.  We then need only translate all the charts of $M$ by an arbitrarily small amount to ensure that all singularities project to rational coordinates; the surface is now square-tiled with a vertical single-cylinder direction.
 
Periodic directions correspond to rational directions (now that the singularities all project to rational coordinates), so the flow in any irrational direction will be completely aperiodic.  By choosing a direction arbitrarily close to vertical, but still irrational, and labeling this direction as $\alpha_m$, we obtain our countable, dense set.  See \autoref{figure - dense collection} to see the effect of this minor shear in direction.

In the language of zippered rectangles, we may make arbitrarily small changes to the lengths of the intervals of the base so that the interval exchange is purely periodic, and without loss of generality we may assume that all lengths are rational.  By also forcing the heights of the rectangles to be rational, the surface is square-tiled.  While this change forbids the surface from the space of zippered rectangles in the vertical direction, which now has saddle connections, an arbitrarily small change to the direction of the flow will remove all saddle connections, so {\em some} representation of this surface appears in the space of zippered rectangles.
\end{proof}
\end{lemma}

We will have no further explicit use for the directions $\alpha_m$; their inclusion is a formality to ensure that our surfaces are proper zippered rectangles.  We will generally be more concerned with the single-cylinder direction(s) on $M_m$.  In this case we will represent single-cylinder directions as `vertical' in figures for convenience.  As such, we will refer to the zippered rectangle $(M_m,\alpha_m)$ simply as $M_m$.  Let $\{M_m\}$ be the countable, dense set of surfaces guaranteed by \autoref{corollary - silly}.

Comparing the surfaces in \autoref{figure - almost single cylinder}, the new surface is nearby within the same stratum.  The flow in direction $\theta$ was periodic for the original surface, and a slight adjustment to the direction of the flow in the new surface will yield a cylinder which is of large measure, periodic in a direction very close to $\theta$, with nearly the same period.  These claims (direct from \autoref{figure - almost single cylinder}) are summarized for convenience in the following lemma.

\begin{lemma}\label{lemma - small balls}
Let $M_{m_0}$ be some particular surface taken from the countable, dense subset given by \autoref{corollary - silly}, with a single-cylinder direction given by $\theta$.  Then for any $\epsilon>0$, there is some $\delta>0$ (depending not just on $\epsilon$ but also on the choice of $\theta$) such that for any $(M', \alpha') \in B_{\delta}(M_{m_0},\alpha_{m_0})$, we have {\em some} direction $\theta^*$ on $M'$ for which there is a set $A \subset M'$ with the following properties:
\begin{itemize}
\item $A$ is a periodic cylinder in direction $\theta^*$, of period $h^*$ and width (measured perpendicular to $\theta^*$) $w^*$,
\item $\mu(A) > (1- \epsilon)\mu(M')$,
\item $|h-h^*| < \epsilon \min\{h,h^*\}$ and $|w-w^*| < \epsilon \min\{w,w^*\}$,
\item $\max\{h,h^*\}| \theta - \theta^*| < \epsilon \min\{w,w^*\}$.
\end{itemize}
\end{lemma}

That is, the surface $(M', \alpha') \in \mathcal{H}(d_1,\ldots,d_n)$ has a direction $\theta'$ which is very nearly a single-cylinder direction: the representation of the flow in direction $\theta'$ is closely approximated by the single-cylinder direction $\theta$ on the surface $M_{m_0}$.  

We now are ready to state our main result:
\begin{theorem}[Residual Generic Ergodicity]\label{theorem - main}
For each $\mathcal{H}(d_1,\ldots,d_n)$, $\mathcal{E}$ is residual within $\mathcal{H}(d_1,\ldots,d_n)$: there is a dense $G_{\delta}$ set of surfaces for which the skew product formed with almost-every collection of cuts is ergodic in almost-every direction.
\end{theorem}

\section{Proof of Residual Generic Ergodicity}
\begin{lemma}
\label{lemma - residual sets}
Let $\{\mathcal{P},\mu\}$ be a $\sigma$-finite measure space, and let $Y$ be a separable metric Baire space such that $\{y_i\}$ $(i=1,2,\ldots)$ is dense in $Y$.  Assume that there is a measurable function
\[\epsilon(i,k,\psi): \mathbb{Z}^+ \times \mathbb{Z}^+ \times \mathcal{P} \longrightarrow \mathbb{R}^+.\]
Define the sets 
\[A_{\psi} = \bigcap_{k=1}^{\infty} \bigcup_{i=1}^{\infty} B_{\epsilon(i,k,\psi)}(y_i),\]
where $B_\epsilon(y) \subset Y$ is the open ball of radius $\epsilon$ centered at $y$.   Then the set 
\[A = \left\{ y \in Y : \mu\left( \{\psi \in \mathcal{P} : y \notin A_{\psi} \} \right) =0\right\}\] is residual.
\begin{proof}
Note that the $A_{\psi}$ are dense $G_{\delta}$ sets.  Without loss of generality, we may assume that $\mu(\mathcal{P})=1$; if the result holds in this circumstance, then the assumption that $\{\mathcal{P}, \mu\}$ is $\sigma$-finite and the fact that a countable intersection of residual sets in a Baire space is residual completes the proof.  It then suffices to prove that for all $\delta>0$, the set
\[A(\delta) = \left\{ y \in Y : \mu \left( \{\psi \in \mathcal{P} : y \in A_{\psi} \} \right) \geq 1-\delta\right\}\] is residual; if this is shown, then we may intersect $A(1/n)$ for $n=1,2,\ldots$.

So let $\delta>0$.  For each $i$ and $k$, define $\epsilon^*(i,k)>0$ so that
\[ \mu \left( \{\psi \in \mathcal{P} : \epsilon(i,k,\psi) < \epsilon^*(i,k) \} \right) \leq \delta \cdot 2^{-k}.\]
Then the set
\[A^* = \bigcap_{k=1}^{\infty} \bigcup_{i=1}^{\infty} B_{\epsilon^*(i,k)}(y_i)\] is a dense $G_{\delta}$ set.  Let $y \in A^*$.  Then for each $k$, there is some $i$ so that $y \in B_{\epsilon^*(i,k)}(y_i)$, which implies that
\[ \mu \left( \{\psi \in \mathcal{P} : y \notin B_{\epsilon(i,k,\psi)}(y_i) \}\right) \leq \delta \cdot 2^{-k}.\]
By intersecting the above over $k=1,2,\ldots$, we have $\mu(A^*) \geq 1-\delta$, or $A^* \subseteq A(\delta)$.
\end{proof}
\end{lemma}
That is, if for almost every parameter $\psi$ the set $A_{\psi}$ is residual, then there is a residual set of $y$ which belong to almost-every $A_{\psi}$.  It is quite important that the set $\{y_i\}$ (the centers of the balls) be invariant with respect to different $\psi \in \mathcal{P}$; otherwise, the result may not hold.  For example, one may readily construct second-countable $Y$ and $\{\mathcal{P}, \mu\}$ a probability space so that for almost every $\psi \in \mathcal{P}$ there is some associated residual set $A_{\psi} \subset Y$, but for \emph{every} $y \in Y$, $\mu \left( \psi \in \mathcal{P} : y \in A_{\psi} \right) =0$.

Our $\sigma$-finite parameter space in applying \autoref{lemma - residual sets} will be 
\[ \mathcal{P} = \mathbb{S}^1 \times \Gamma,\]
where $\Gamma$ was given in \eqref{eqn - space of cuts}.  A point in this space corresponds to a choice of direction $\theta$ for geodesic flow $\tilde{\varphi}_t$, combined with a collection of cuts $\gamma\in \Gamma$.  Our Baire space $Y$ will be a stratum of translation surfaces (\autoref{lemma - baire}).

Let the stratum $\mathcal{H}(d_1, \ldots, d_n)$ be fixed, and denote by $\{M_m\}$ the countable, dense collection of square-tiled translation surfaces which admit some single-cylinder direction from \autoref{corollary - silly}.

\begin{lemma}\label{lemma - single cylinder directions}
For almost-every direction $\theta$, for each $M_m$ there is a sequence $\{\theta_n\}$ of single-cylinder directions of period $h_n$ and width $w_n$ on $M_m$ such that
\begin{equation}\label{eqn - good approx} \lim_{n \rightarrow \infty} \frac{h_n}{w_n}\|\theta_n- \theta\| =0. \end{equation}
\begin{proof}
As there are only countably many $M_m$, it suffices to prove the claim for a single $M_m$, which follows from \cite[Lemma 2.6]{ergodic-infinite}.
\end{proof}
\end{lemma}

\begin{figure}[th]
\center{\begin{tikzpicture}
\draw (0,0) rectangle (2,4);
\filldraw (0,0)[fill=shade2] rectangle (.2,4);
\filldraw (1.9,0)[fill=shade2] rectangle (2,4);
\draw [xshift=-10pt](0,.5) edge[->] node[midway, xshift=-20pt]{$\theta_{n_j}$} (0,3.5) ;
\draw [dotted, ->] (1.05,0) --  node{$\cdot$} node[right]{$d$} (1.05,1.5)--(1.05,4);
\node at (.6,2)[below]{$L_j$};
\node at (1.5,2)[below]{$R_j$};
\draw [decorate, decoration={brace,amplitude=10pt, mirror},xshift=0pt,yshift=0pt] (2.2,.2)-- (2.2,3.8) node[midway, xshift=22pt]{$h_{n_j}$};
\end{tikzpicture}}
\caption{\label{figure - avoiding} By assumption the only element of the set $D$ not within the gray region is the point $d$, which may then be placed in the center of a sub-cylinder of large measure which avoids all singularities and all other elements of $D$.  As $j \rightarrow \infty$ the grayed-out region tends to zero measure, and thus $d$ tends towards the center of the cylinder.  The remaining sub-cylinders are called orbit neighborhoods of $d$, with left and right halves $L_j$ and $R_j$.}
\end{figure}
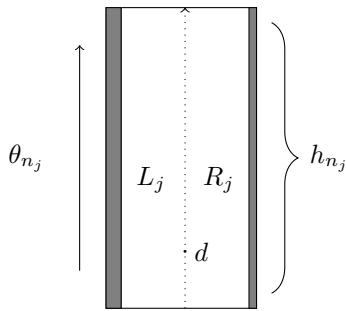

We recall \cite[Def. 3.1]{ergodic-infinite} and strengthen:

\begin{definition}\label{definition - avoiding set}
Let $M$ be a fixed compact translation surface which has a dense set of single-cylinder directions, let $\theta$ be some direction, and let $D$ be a finite subset of $M$.  We will say that the set $D$ is \emph{uniformly self-avoiding} in direction $\theta$ if there is a sequence of single-cylinder directions $\theta_n \rightarrow \theta$ such that for each $d \in D$ there is a subsequence $\theta_{n_j}$, along which we may place $d$ at the center of an invariant subset of the single cylinder, the complement of which tends towards zero measure, which contains no singularities of $M$ and no other points of $D$.  See \autoref{figure - avoiding}.  The two sides of this set, separated by the orbit of $d$, are called the \emph{left} and \emph{right} halves of the \emph{orbit neighborhoods} of $d$.
\end{definition}

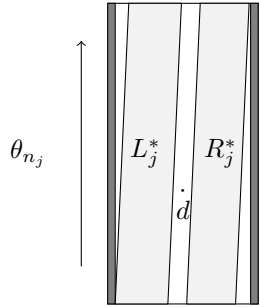
\begin{figure}[bh]
\center{\begin{tikzpicture}
\draw (0,0) rectangle (2,4);
\filldraw (0,0)[fill=shade2] rectangle (.1,4);
\filldraw (1.9,0)[fill=shade2] rectangle (2,4);
\draw [fill=shade4, xslant=.045] (.1,0) rectangle (.8,4);
\draw [fill=shade4, xslant=.045] (1.05,0) rectangle (1.7,4);
\draw [xshift=-10pt](0,.5) edge[->] node[midway, xshift=-20pt]{$\theta_{n_j}$} (0,3.5) ;
\node at (1,1.5) {$\cdot$};
\node at (1,1.5) [below]{$d$};
\node at (.5,2){$L^*_j$};
\node at (1.5,2){$R^*_j$};
\end{tikzpicture}}
\caption{\label{figure - neighborhoods}The flow in direction $\theta$ is slightly skewed compared to direction $\theta_{n_j}$, but sufficiently close to allow for the construction of the sets $L^*_j$ and $R^*_j$, each of measure very nearly one half the total measure.}
\end{figure}

\begin{lemma}\label{lemma - avoiding}
For any $N \in \mathbb{N}$, almost every choice of $N$ points in the plane lifts to every $M_m$ as a uniformly self-avoiding set.  Furthermore, if $M_{m_0}$ is fixed, and $\theta$ is an irrational direction which admits a sequence $\{\theta_n\}$ of single-cylinder directions (with periods $h_n$) which satisfies \autoref{eqn - good approx}, then for each $d$ in this self-avoiding set, if we pass to the subsequence $n_j$ from \autoref{definition - avoiding set} and denote the left and right halves of a sequence of orbit neighborhoods corresponding to $d$ by $L_j$ and $R_j$, then we have subsets $L^*_j \subset L_j$, $R^*_j \subset R_j$ such that as $j \rightarrow \infty$,
\begin{itemize}
\item $\mu(L^*_j) \rightarrow (1/2)\mu(M_{m_0})$,
\item $\sup_{x \in L^*_j}\|\varphi^{\theta}_{h_{n_j}}(x) - x\| \rightarrow 0$,
\item for any fixed $t_0 \in \mathbb{R}^+$, $\mu \left( \varphi^{\theta}_{t_0} L^*_j \triangle L^*_j \right) \rightarrow 0$,
\end{itemize}
and similarly for $R^*_j$. 
\begin{proof}
Note that the flows are in direction $\theta$, not $\theta_{n_j}$.  For the first claim, as there are countably many $M_m$, it suffices to prove that almost every choice of $N$ points lifts to any particular $M_{m_0}$ as a uniformly self-avoiding set.  We rescale a fixed $M_{m_0}$ to be a $k$-fold cover of the unit square and denote the affine image of $\theta_n$ by $p_n/q_n$ (periodic directions on covers of the unit square must be rational).  Single cylinders therefore have width $1/q_n$ and height $kq_n$ for some fixed $k \in \mathbb{Z}^+$.  The first claim now follows from \cite[Thm 4.1]{MR0419394}: for \emph{any} sequence $q_n \rightarrow \infty$, for almost every $x$ we have the sequence $kq_n x$ uniformly distributed modulo one.  So if there are $N$ points $d_1,d_2,\ldots,d_N$, let $\pi_x(d_i)$ be the $x$-coordinate of the projection to the plane.  So for each $i' \in\{1,2,\ldots,n\}$ we can create a subsequence $p_{n_j}/q_{n_j}$ such that modulo 1 $q_{n_j}\pi_x(d_{i'})$ tends towards $1/2$ as $j \rightarrow \infty$ (and therefore $d_{i'}$ is arbitrarily close to the center of the cylinder) while modulo 1 every other $i \neq i'$, $q_{n_j}\pi_x(d_{i})$ tends to either zero or one (and therefore all other $d_{i}$ are arbitrarily close to the boundary of the cylinder).

Now let $d$ be an element of a uniformly self-avoiding set.  We form the sets $L^*_j$ and $R^*_j$ by removing a small neighborhood of the orbit of $d$ (in the direction $\theta$).  See \autoref{figure - neighborhoods}.  Using \autoref{eqn - good approx} and \autoref{definition - avoiding set}, we see that $L^*_j$ and $R^*_j$ will tend to half of the total measure.  That the sets are nearly invariant under the flow in direction $\theta$ (pointwise through time $h_{n_j}$ and as a set through any fixed time $t_0$) follows from \autoref{eqn - good approx} and the fact that $h_{n_j}$ is diverging, respectively: the flow for any fixed $t_0$ will produce a negligible change in either $L^*_j$ or $R^*_j$.
\end{proof}
\end{lemma}

Compare to \cite[Cor. 3.4]{ergodic-infinite}.  In that result only a specific square-tiled surface was considered, so any positive measure for the sets $L$, $R$ sufficed.  We have strengthened that approach to provide \emph{the same positive limiting measure} for the sets $L^*_j$ and $R^*_j$ for every $M_m$.

\begin{corollary}\label{corollary - good subset of cuts}
For almost every collection of cuts $\gamma \in \Gamma$, the lift of $\gamma$ to each $M_m$ is unique, and for almost every direction $\theta$ the {\em starting} points $\textbf{x}_i^{\pm}$ of the segments of the lift of $\gamma$ are a uniformly self-avoiding set in the direction $\theta$ on every $M_m$.

\begin{proof}
Uniqueness of the lift of $\gamma$ to each $M_m$ merely requires that the geodesic from $p'$, the projection of the marked point $p \in M_{m_0}$, to each $(x_i^{\pm},y_i^{\pm})$, as well as the line segments originating from these points, not encounter a singularity when lifted to $M$.  As there set of singularities is of zero measure (it is finite), generic cuts lift uniquely to \emph{each} $M_m$.  

It follows from \cite{Patterson, MR688349} (see also \cite[Lemma 2.6]{ergodic-infinite}) that for any fixed $M_{m_0}$, almost every $\theta$ admits a sequence $\{\theta_n\}$ as in \autoref{lemma - avoiding}, so for each $M_m$, for almost every $\theta$ we have that almost every collection of cuts lifts to $M_m$ uniquely (as outlined above) and with endpoints which are a uniformly self-avoiding set in direction $\theta$.  Fubini's theorem allows us to rephrase: almost every collection of cuts lifts uniquely as a self-avoiding set to any fixed $M_{m_0}$ in almost-every direction.  Since the collection $\{M_m\}$ is countable, it follows that almost every collection of cuts lifts to a uniformly self-avoiding set on every $M_m$ in almost every direction.
\end{proof}
\end{corollary}

\begin{definition}
We call the pair $(\gamma,\theta)$ {\em generic} if it satisfies the conclusion of Corollary \ref{corollary - good subset of cuts}.
\end{definition}

To summarize our progress so far: we have a countable, dense collection of surfaces $M_m$ in each stratum.  For almost-every direction $\theta$, each one of these surfaces has a sequence of single-cylinder directions converging in a very strong manner to $\theta$.  For almost-every collection of cuts, for each starting point of each cut, on each $M_m$ we may pass to a subsequence of these directions which isolates that starting point from the orbits of every other starting point in a uniform manner.  We will now show how this isolation of starting points of the cuts leads to ergodic skew products.

An \emph{essential value} of the skew product $\{\tilde{M}, \tilde{\mu}, \tilde{\varphi}_t\}$ is a $g \in G$ such that for any set $A \subset M$ of positive measure and $\epsilon>0$ there is a $t>0$ such that (recall \eqref{eqn - balls for essential values})
\[ \mu \left( A \cap \varphi_t(A) \cap S_t^{-1}(B_{\epsilon}(g))\right) >0.\]
The set of essential values forms a closed subgroup of $G$, and the skew product is ergodic if and only if the flow $\{M, \mu, \varphi_t\}$ is ergodic and the set of essential values is all of $G$ (see e.g. \cite{A}).  Clearly, for $\{ \tilde{M}, \tilde{\mu}, \tilde{\varphi}_t\}$ to be ergodic it is necessary that the values $f$ takes on the cuts to generate a dense subgroup of $G$.  Our technique for finding essential values is the following:

\begin{lemma}\label{lemma - big list}
Suppose that $M$ is a translation surface with geodesic flow $\{M, \mu, \varphi_t\}$ in some direction $\theta$, with a collection of cuts $\gamma$ and corresponding skew product $\{\tilde{M}, \tilde{\mu}, \tilde{\varphi}_t\}$.  For $j=1,2,\ldots$, assume we have times $t_j \rightarrow \infty$ and sets $L^*_j$ such that

\begin{itemize}
\item $L^*_j$ contains no singularities of $M$ and no endpoints of any cut,
\item $\lim \mu(L^*_j) = (1/2)\mu(M)$,
\item $\sup_{x \in L^*_j} \left\|\varphi^{\theta}_{t_n}(x) - x\right\| \rightarrow 0$,
\item for any fixed $t_0 \in \mathbb{R}^+$, $\mu \left( L^*_j \triangle \varphi^{\theta}_{t_0} L^*_j \right) \rightarrow 0$, and
\item for all $x \in L^*_j$, $S_{t_j}(x)=f(\gamma_i)$.
\end{itemize}
Then the skewing values $\pm g_i = f(\gamma_i^{\pm})$ taken on the cuts $\gamma_i^{\pm}$  are essential values of $\{\tilde{M}, \tilde{\mu}, \tilde{\varphi}_t\}$.
\begin{proof}
Note that we have required the ergodic sums $S_{t_j}(x)$ to be constant on each $L^*_j$ and $R^*_j$, which each form a sequence of \emph{quasi-rigidity sets} \cite[Def. 2.1]{ergodic-infinite}, and we directly have each of $\pm g_i$ as essential values by applying \cite[Lemma 2.2]{ergodic-infinite}.
\end{proof}
\end{lemma}

Compare to \autoref{lemma - avoiding} and \autoref{figure - neighborhoods}: the last condition in \autoref{lemma - big list} is the last remaining condition to verify in order to begin producing essential values for our skew products.

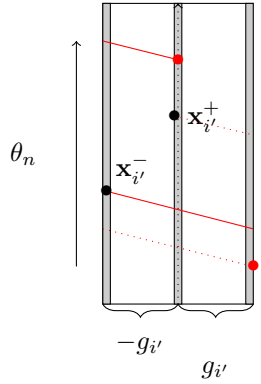
\begin{figure}[tbh]
\center{\begin{tikzpicture}
\draw (0,0) rectangle (2,4);
\filldraw (0,0)[fill=shade3] rectangle (.1,4);
\filldraw (1.9,0)[fill=shade3] rectangle (2,4);
\draw [xshift=-10pt](0,.5) edge[->] node[midway, xshift=-20pt]{$\theta_{n}$} (0,3.5) ;
\filldraw (.95,0)[fill=shade3] rectangle (1.05,4);
\draw [dotted, ->] (1,0) -- (1,4);
\node at (.05,1.8)[right]{$\textbf{x}^-_{i'}$};
\node at (1,2.5)[right]{$\textbf{x}^+_{i'}$};
\draw [color=red] (.05,1.5) -- (2,1);
\draw [color=red] (0,3.5) -- (1,3.25);
\draw [color=red, dotted] (.95, 2.5) -- (2,2.25);
\draw [color = red, dotted] (0,1) -- (2,.5);
\foreach \Point in {(.05,1.5),(.95,2.5)}
    \node at \Point {\textbullet};
\foreach \Point in {(1,3.25),(2,.5)}
    \node [color=red] at \Point {\textbullet};
\draw [decorate, decoration={brace, amplitude=5pt, mirror}] (0,0)--(1,0) node[midway, yshift=-15pt]{$-g_{i'}$};
\draw [decorate, decoration={brace, amplitude=5pt, mirror}] (1,0)--(2,0) node[midway, yshift=-25pt]{$g_{i'}$};
\end{tikzpicture}}
\caption{\label{figure - nice sums} For generic cuts and generic directions, on each particular $M_{m_0}$ we may pass to a subsequence of single-cylinder directions on which the cuts traverse the single cylinder in a regular manner.  All cuts other than a chosen $\gamma^{\pm}_{i'}$ have both endpoints restricted near the boundary and therefore contribute nothing to the ergodic sums in the cylinder, while the chosen cut $\gamma^{\pm}_{i'}$ splits the cylinder into two halves with very precise ergodic sums.}
\end{figure}

\begin{lemma}\label{lemma - uniform bound}
Fix generic $(\gamma , \theta) \in \Gamma \times \mathbb{S}^1$, and let $f:\gamma \rightarrow G$, where $G$ is a locally compact Abelian group.   Then for each $M_{m_0}$ we have that for each $\gamma^{\pm}_i$ with value $f(\gamma_i^{\pm})$ we may pass to a subsequence $\{\theta_{n_j}\}$ of single cylinder directions so that we may find a sequence of sets $L^*_j$ to satisfy \autoref{lemma - big list}.
\begin{proof}
The proof is a slightly more refined application of the same use of \cite[Thm. 4.1]{MR0419394} that we used in \autoref{lemma - avoiding}.  Pick a generic $(\gamma,\theta)$, 
and let $M_{m_0}$ be fixed.  From \autoref{lemma - single cylinder directions} we obtain a sequence $\theta_n$ of single-cylinder directions on $M_{m_0}$ converging to $\theta$.  For each pair of cuts $\gamma^{\pm}_i$, in direction $\zeta_i$ and 
with length $\ell_i$, project the length of $\gamma_i$ in the direction $\theta$ and denote this projected length $\pi(\ell_i)$.  Let $w_{n}$ be the width of the single-cylinder on $M_{m_0}$ in direction $\theta_n$.  

Pick a particular pair of cuts $\gamma_{i'}^{\pm}$.  From \autoref{lemma - avoiding}, let us assume without loss of generality that our sequence $\theta_n$ already isolates the starting point $\textbf{x}_{i'}^+$ of $\gamma_{i'}^+$ from the starting point of every other cut (including $\textbf{x}_{i'}^-$).  For each $i \neq i'$ (and $\textbf{x}_{i'}^-$), we pass to successive subsequences to guarantee that the {\em endpoints} of  the cuts $\gamma_i^{\pm}$  are also forced to be  near the boundary of the cylinders.  More precisely we require that  $(\pi(\ell_i)/w_{n_j}) \bmod 1$  converges to zero/one as $j \rightarrow \infty$, i.e.\ the cuts $\gamma^{\pm}_i$ cross the cylinders very nearly an exact integer number of times.  As the projected length $\pi(\ell_i)$ is equal for both $\gamma_i^+$ and $\gamma_i^-$ this subsequence automatically isolates the endpoint of $\gamma_i^-$ near the boundary as well.  Since there are only finitely many $i \neq i'$, we may presume, then, that {\em every} cut other than $\gamma_{i'}$ has this property: both endpoints of each $\gamma_i^{\pm}$ are trapped near the boundaries of the cylinders.

So now consider the chosen pair of cuts $\gamma_{i'}^{\pm}$, with starting points $\textbf{x}_{i'}^{\pm}$.  Since we began with a subsequence of directions according to \autoref{definition - avoiding set}, we have that $\textbf{x}_{i'}^+$ is placed near the centers of the cylinders, and $\textbf{x}_{i'}^-$ is trapped near the boundary along with all the other starting points.  Since the length $\ell_{i'}$ of this particular pair of cuts was chosen generically, by passing to a subsequence once more we may force the {\em endpoint} $\varphi_{\ell_{i'}}^{\zeta_{i'}}(\textbf{x}_{i'}^+)$ to be trapped near the {\em boundary} of the cylinder.  Equivalently, we may generically force $(\pi(\ell_{i'})/w_{n_j}) \bmod 1$ to converge to $1/2$: the length of this cut traverses the cylinders very nearly an integer number of times, {\em plus one half}.  

So for $i \neq i'$, the cuts begin very near the boundary of the single cylinders and are of length (projected to direction $\theta$, which is very near to the directions $\theta_{n_j}$) very nearly an integer times the width of the cylinder.  Away from the boundary of the cylinder, then, \emph{the ergodic sums of any point flowing the height of the cylinder see no contribution from the cut $\gamma^{\pm}_i$}!  As the pair of cuts $\gamma^{\pm}_i$ are of equal length, they must traverse the cylinder an equal number of times.  Since we have $f(\gamma^-_i)=-f(\gamma^+_i)$, the net effect is that these cuts do not contribute to the ergodic sums for most points in this direction through the height of the cylinder.

For the remaining cut $\gamma^{\pm}_{i'}$, the endpoint of $\gamma^+_{i'}$ is very near the boundary of the cylinder, while the endpoint of $\gamma^-_{i'}$ is very near the \emph{middle} of the cylinder (a reversal of the placement of the starting points $\textbf{x}^{\pm}_{i'}$, as the projected length $\pi(\ell_{i'})$ is very nearly an integer-plus-one-half times the width of the cylinder).  So we find that the ergodic sums neatly split the cylinder; away from the boundary and away from the middle, the only two values seen are exactly $g_{i'}=f(\gamma_{i'})$ on the right half of the cylinder and $-g_{i'}=-f(\gamma_{i'})$ on the left.  See \autoref{figure - nice sums}; the sets $L^*_j$ and $R^*_j$ are obtained from this representation in the same way as in \autoref{figure - neighborhoods}.
\end{proof}
\end{lemma}

We finally turn to the question of the appropriate $\epsilon$ to use as the radius of the balls in constructing our $G_{\delta}$ sets to complete the proof of \autoref{theorem - main}.  Let $\psi = (\theta,\gamma) \in \mathcal{P}$ be generic in the sense that all generic results proved up to this point hold simultaneously (\autoref{lemma - single cylinder directions}, \autoref{lemma - avoiding}, \autoref{corollary - good subset of cuts}, \autoref{lemma - big list}, and \autoref{lemma - uniform bound}).

Let $k$ be a fixed positive integer, and fix some $M_{m_0}$.  For each $d$, an endpoint of a cut, we have a sequence of single-cylinder directions $\theta_{n}$ according to \autoref{lemma - avoiding}.  We may therefore find {\em some} single-cylinder direction $\theta_{n_j}$ on $M_{m_0}$ which produces two orbit-neighborhoods $L_j^*$ and $R_j^*$ according to \autoref{figure - neighborhoods} such that each has measure at least $(k-1)/2k$ and each is of height at least $k$.  Fix this choice of single-cylinder direction $\theta_{n_j}$ (which depends not only on $m_0$ and $\theta$, but also $k$).

We let $\epsilon(m_0,k,d)>0$ be small enough (recall \autoref{lemma - small balls}, making the size of the balls in the stratum smaller if necessary) that any surface within $\epsilon(m_0,k)$ of each $(M_{m_0},\alpha_{m_0})$ satisfies:
\begin{enumerate}
\item the lift of $\gamma$ to all nearby $M$ in this neighborhood will be unique,
\item excluding a set of measure at most $\mu(M)/k$, there is a periodic cylinder in direction $\theta'$ of period $h'$ such that $\max{(h,h')}\cdot|\theta_{n_j}-\theta'| < w/(2k)$ for all points in this cylinder, where $w$ is the width of the cylinder.  See \autoref{figure - almost single cylinder}, right.  
\item for those $x \in M$ of within this large set and also not within distance $w/(2k)$ of any lift of a cut, we have
\[S_{h'}(x,M,\theta) = S_{h'}(x,M,\theta')=S_h(\hat{x},M_{m_0},\theta_{n_j})=S_h(\hat{x},M_{m_0},\theta),\]
\end{enumerate}
where $\hat{x}$ denotes the lift to $M_{m_0}$ of the projection of $x \in M$ to the plane, and $\theta_{n_j}$ is a single-cylinder direction on $M_{m_0}$ very near $\theta$.  Informally, point-by-point these conditions amount to:
\begin{enumerate}
\item varying the edges of the polygon will not cause any cut to cross a singularity,
\item the flows in the two directions $\theta$ and $\theta'$ through times $h$ or $h'$ are more-or-less indistinguishable for most points, and also virtually indistinguishable from the single-cylinder direction $\theta_{n_j}$, and
\item for most of the surface $M$, the ergodic sum through time $h'$ in direction $\theta$ is identical to the ergodic sum through time $h$ in direction $\theta$ of the corresponding point on $M_{m_0}$.
\end{enumerate}

The equality in the last item is accomplished by starting with $M_{m_0}$ in direction $\theta$, approximating by a single-cylinder direction $\theta_{n_j}$, choose a surface in the $\epsilon$-ball around $M_{m_0}$ by making small adjustments to the vectors which define the saddle connection in direction $\theta_{n_j}$ to form the surface $M$, finding a cylinder of large measure on the resulting surface in direction $\theta'$, and then returning the flow to direction $\theta$.  At each stage, ergodic sums are maintained pointwise except on sets of arbitrarily small measure, according to the diameter of the ball used.

Having chosen the cuts $\gamma$ and the direction of the flow $\theta$, we therefore have a small ball around $M_{m_0}$ for which any surface in that ball is very closely approximated by the single-cylinder direction $\theta_{n_j}$ on $M_{m_0}$, at least through time $h_{n_j} \geq k$.  As $k$ increases, the strength of this approximation necessarily improves in two ways: increasing height of the cylinders forces $\theta_{n_j} \rightarrow \theta$, and the measure of the orbit neighborhoods which directly transfer to any $M$ in the ball is tending towards $1/2$.  So with $m_0$ fixed, for each $k$ and each endpoint of each cut $d$, we obtain an allowable radius $\epsilon(m_0,k,d)$.  Define $\epsilon(m_0,k, \psi)$ to be the minimum of these radii across the finite number of endpoints of the cuts (recall the original dependence on $\psi$, the collection of cuts and the direction of geodesic flow).

Let
\[A_{\psi} = (E(\theta)) \cap \left(\bigcap_{k\in \mathbb{N}} \bigcup_{i \in \mathbb{N}} B_{\epsilon(i,k,\psi)}(M_m,\alpha_m)\right),\]
where $E(\theta)$ is the set of $M$ for which the flow in direction $\theta$ is ergodic, and $\theta$ was the direction of geodesic flow specified in $\psi$.
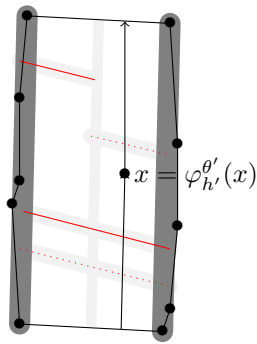
\begin{figure}[tbh]
\center{\begin{tikzpicture}
\draw [line width=4, color=shade4] (.95,-.04)--(1.05,4.05);
\draw [line width=6, color=shade4, cap=round] (.95,2.5)--(2,2.25);
\draw [line width=6, color=shade4, cap=round] (0,1)--(2,.5);
\draw [line width=6, color=shade4, cap=round] (.05,1.5)--(2,1);
\draw [line width=6, color=shade4, cap=round] (0,3.5)--(1,3.25);
\draw [line width=8, color=shade2, cap=round] (1.9,-.1)--(2,4);
\draw [line width=8, color=shade2, cap=round] (.1,4.1)--(0,0);
\filldraw[color=shade2](0,0)--(-.1,1.6)--(0,1.9)--(0,3)--(.1,4.1)--(0,0);
\filldraw[color=shade2](2,4)--(2.1,2.4)--(2.1,1.3)--(2,.2)--(1.9,-.1)--(2,4);
\draw [color=red] (.05,1.5) -- (2,1);
\draw [color=red] (0,3.5) -- (1,3.25);
\draw [color=red, dotted] (.95, 2.5) -- (2,2.25);
\draw [color = red, dotted] (0,1) -- (2,.5);

\node at (1.4,2){\textbullet};
\node at (1.4,2)[right]{$x=\varphi^{\theta'}_{h'}(x)$};
\draw [->] (1.4,2)--(1.4,4.025);
\draw [->] (1.35,-.07)--(1.4,2);
\draw (0,0)--(-.1,1.6)--(0,1.9)--(0,3)--(.1,4.1)--(2,4)--(2.1,2.4)--(2.1,1.3)--(2,.2)--(1.9,-.1)--(0,0);
\foreach \Point in {(0,0), (-.1,1.6), (0,1.9), (0,3),(.1,4.1),(2,4), (2,.2), (2.1,1.3), (2.1,2.4), (1.9,-.1)}
    \node at \Point {\textbullet};
\end{tikzpicture}}
\caption{\label{figure - make epsilon} If $M$ is sufficiently close to $M_{m_0}$, then `most' of $M$ looks like a single cylinder (all but the dark-gray region, as in \autoref{figure - almost single cylinder} (right) and \autoref{figure - avoiding}), and the cuts are in almost the same place, with the same number of crossings, so that for `most' of this periodic cylinder (now also excluding the light-gray regions) the ergodic sums are identical to corresponding points on $M_{m_0}$ (as in \autoref{figure - nice sums}).  
}
\end{figure}
\begin{lemma}
Let $M \in A_{\psi}$, and let $\tilde{M} = M \times G$ with geodesic flow in direction $\theta$ defined as the skew product over $M$ given by \eqref{eqn - sums}.  Then the geodesic flow on $\tilde{M}$ in direction $\theta$ is ergodic.
\begin{proof}
The set $E(\theta)$ is residual in each stratum for any $\theta$ \cite[Corollary 2]{1986}, so $A_{\psi}$ is residual as well.  For $M \in A_{\psi}$, then, the flow in direction $\theta$ is ergodic.  It remains to show that each value $g \in G$ which is the value assigned to some cut is an essential value of the flow $\{\tilde{M}, \tilde{\mu}, \tilde{\varphi}_t\}$; since these values generate a dense subgroup of $G$, this finding will complete the proof.

Pick an arbitrary $g \in \{g_1,\ldots,g_N\}$ (the values taken on the cuts), and let $\textbf{x}_i^+$ be the starting point of a cut which takes value $g$.  Create a sequence of approximations as in \autoref{figure - make epsilon}.  By construction, as $M \in A_{\psi}$, we have a sequence of single-cylinder surfaces $M_m$ with corresponding representations as in \autoref{figure - make epsilon}, from which we may pull back to $M$ a set of measure arbitrarily close to $1/2$ on which we have ergodic sums constantly equal to $g$; via \autoref{lemma - uniform bound} the ergodic sums on the corresponding sets in $M_m$ are constant, and our balls were picked small enough that most of these sets transfer to $M$ without affecting the ergodic sums.  The value $g$, then, is an essential value of the flow in direction $\theta$ on $M$ (\autoref{lemma - big list}).  Since the choice of $g$ was arbitrary and we assumed that $\{g_1,\ldots,g_N\}$ generate a dense subgroup of $G$, it follows that the flow on $\tilde{M}$ in direction $\theta$ is ergodic.
\end{proof}
\end{lemma}

So for any generic choice of direction $\theta$ and cuts $\gamma$, we form a residual set of surfaces $M$ for which the skew product transformation on $\tilde{M}$ is ergodic in direction $\theta$.  By applying \autoref{lemma - residual sets}, this completes the proof of \autoref{theorem - main}.

\bibliography{residual-bibfile2}
\bibliographystyle{plain}
\end{document}